\documentclass[letterpaper, 10 pt, conference]{ieeeconf}  

\IEEEoverridecommandlockouts                              

\usepackage{lipsum}
\usepackage[noadjust]{cite}  
\usepackage{amsmath,amssymb}
\usepackage{aligned-overset}
\usepackage{enumerate}
\usepackage{bm}
\usepackage{amsfonts}
\usepackage{anyfontsize}
\usepackage{cancel}

\usepackage{float}

\usepackage{multirow}
\usepackage{array}

\newtheorem{theorem}{Theorem}[section]
\newtheorem{lemma}{Lemma}[section]
\newtheorem{proposition}{Proposition}[section]

\newtheorem{remark}{Remark}[section]
\newtheorem{assumption}{Assumption}[section]
\newtheorem{problem}{Problem}[section]

\newcommand\numberthis{\addtocounter{equation}{1}\tag{\theequation}}
\makeatletter

\usepackage{algpseudocode}
\usepackage{algorithm}

\usepackage{setspace}

\allowdisplaybreaks  

\usepackage[unicode]{hyperref}
\hypersetup{
    colorlinks,
    citecolor=black,
    filecolor=black,
    linkcolor=black,
    urlcolor=black
}

\usepackage{nccmath} 

\usepackage{graphicx}

\usepackage[font=sf]{caption}
\usepackage{subcaption}
\usepackage{mathrsfs}   

\newcommand{\me}{\mathrm{e}}
\newcommand{\md}{\mathrm{d}}

\usepackage{accents}

\usepackage{pgfplots}
\pgfplotsset{compat=1.15}
\usetikzlibrary{patterns}
\usetikzlibrary{bending}
\usetikzlibrary{calc}
\usetikzlibrary{arrows,backgrounds,positioning,fit,petri}
\usetikzlibrary{decorations.pathreplacing,decorations.pathmorphing,decorations.markings,decorations.shapes}
\usetikzlibrary{arrows.meta}
\usetikzlibrary{quotes,angles}
\usetikzlibrary{chains}
\usetikzlibrary{hobby}
\usepackage{tkz-euclide}
\definecolor{ttttff}{rgb}{0.2,0.2,1}
\definecolor{ududff}{rgb}{0.30196078431372547,0.30196078431372547,1}
\definecolor{xfqqff}{rgb}{0.4980392156862745,0,1}
\definecolor{qqzzff}{rgb}{0,0.6,1}
\definecolor{ffxfqq}{rgb}{1,0.4980392156862745,0}
\definecolor{qqqqff}{rgb}{0,0,1}
\definecolor{ffqqff}{rgb}{1,0,1}
\definecolor{qqzzqq}{rgb}{0,0.6,0}
\definecolor{ffqqqq}{rgb}{1,0,0}
\definecolor{qqwuqq}{rgb}{0,0.39215686274509803,0}
\definecolor{zzttqq}{rgb}{0.6,0.2,0}
\definecolor{wqwqwq}{rgb}{0.3764705882352941,0.3764705882352941,0.3764705882352941}
\definecolor{qqttcc}{rgb}{0,0.2,0.8}
\definecolor{windBlue}{rgb}{0.3176470588,0.9254901961,0.9450980392}

\newcommand\norm[1]{\| #1 \|}

\renewcommand{\thetable}{\arabic{table}}

\makeatletter
\newcommand{\pushright}[1]{\ifmeasuring@#1\else\omit\hfill$\displaystyle#1$\fi\ignorespaces}
\newcommand{\pushleft}[1]{\ifmeasuring@#1\else\omit$\displaystyle#1$\hfill\fi\ignorespaces}
\makeatother

\newcommand\scalemath[2]{\scalebox{#1}{\mbox{\ensuremath{\displaystyle #2}}}}

\title{\LARGE \bf Integral Compensation for Discrete-time Bearing-only Target Localization and Circumnavigation}

\author{Donglin Sui$^{1}$ and Mohammad Deghat$^{1}$
\thanks{$^{1}$The authors are with the School of Mechanical and Manufacturing Engineering, University of New South Wales, Sydney NSW 2052, Australia (e-mail: {\tt\small d.sui@unsw.edu.au, m.deghat@unsw.edu.au}).}%
}

\begin{document}

\maketitle
\thispagestyle{empty}
\pagestyle{empty}

\begin{abstract}
    This paper investigates the discrete-time control problem of driving an agent to localize and circumnavigate an unknown target at a prescribed radius using bearing measurements. Discretizing existing continuous-time circumnavigation controllers can introduce steady-state orbit bias, even at sampling rates that preserve stability. To eliminate the sampling-induced orbit bias, we propose an adaptive compensator for a circumnavigation controller formulated directly in discrete time. Driven by the estimated agent-to-target distance, the proposed compensator dynamically learns the geometric correction required to achieve the desired circumnavigation radius. Unlike existing static correction methods, neither the proposed controller nor the compensator update explicitly uses the sampling interval. Exponential convergence is rigorously established using Lyapunov-based analysis. Simulation results illustrate the performance of the proposed algorithms.
\end{abstract}


\section{INTRODUCTION}
\label{sec:Introduction}

Target circumnavigation has attracted sustained research interest owing to its practical applications. The task typically requires a single agent or a team of agents to circle a target at a prescribed radius. Earlier studies commonly treat the target's state as available to the agents, an assumption that may fail in practice when an uncooperative target does not disclose this information. Under such circumstances, navigation control must be combined with target localization, leading to a challenging dual control problem \cite{feldbaum1963Dual, rantzer2023Dual}.

The sensing modality becomes particularly important when unknown targets must be localized without compromising stealth. Acquiring distance measurements typically entails transmitting detectable signals toward the target, potentially revealing the agent's presence. Bearing sensing, in contrast, is generally passive: it uses only signals naturally emitted by the target and is feasible with inexpensive devices such as monocular cameras. These advantages have stimulated growing interest in \textit{Bearing-only Target Localization and Circumnavigation} (\textit{BoTLC}).

An early BoTLC architecture in \cite{deghat2010Target,deghat2014Localization} couples a circumnavigation controller with an adaptive target estimator. The controller guides the agent around the target at a specified radius, while the estimator supplies target position estimates using the measured bearings and the agent's own position expressed in an initial reference frame. Extensions address multiple-target localization \cite{chen2023Cooperative,wang2024Achieving,ji2024Adaptive} and unbiased estimation of a constant-velocity target \cite{sui2024Unbiased}. The framework has also been adapted for 3D workspaces \cite{li2018Localization,ma2023Finitetime,ghommam2025Cooperative} and intermittent measurements \cite{zou2021Distributed}, and reformulated to avoid using the agent's own position \cite{wang2022Target,hu2022BearingOnly}. Further developments cover non-holonomic \cite{deghat2012Target} and double-integrator agents \cite{cao2023Bearingonly}, multi-agent systems \cite{boccia2017Tracking,dou2020Target,chen2022Multicircular,sui2024CollisionFree,jiang2026CommunicationFree}, and irregularly shaped circumnavigation orbits \cite{ma2022Finitetime,yue2026FixedTime,cao2021Safe}.

However, all of the studies discussed above formulate and analyze the problem in continuous time, whereas practical sensing systems most commonly provide bearing measurements at discrete times, for example, through monocular cameras. Under ideal conditions, a motion-capture system such as OptiTrack can provide a sampling rate of up to $240$ Hz \cite{schwind2020Effects}, and a straightforward discretization of a continuous-time BoTLC algorithm can be applied without difficulty \cite{sui2025PredefinedTime}. In field applications, however, practical constraints can substantially reduce the effective measurement rate. For example, a computationally constrained quadcopter can typically acquire monocular camera frames at $60-72$ Hz \cite{walter2019UVDAR,zhang2026Observability}, but target recognition using models such as YOLO can further reduce the effective measurement rate to $5-20$ Hz \cite{schilling2021VisionBased,zheng2025Optimal,xu2022OmniSwarm}. These limitations motivate analyzing the estimator-controller framework directly in discrete time rather than relying on discretization of the continuous-time design.

Among the first studies to address this need, \cite{chen2022DiscreteTime} shows that naive forward Euler or Runge--Kutta discretization of continuous-time algorithms leads to divergence at low sampling frequencies. Even at higher control frequencies for which the discretized system remains stable, the uncorrected discretized controller exhibits a sampling-induced steady-state error (SSE). To eliminate this error, \cite{chen2022DiscreteTime} introduces a static geometric correction that depends explicitly on the sampling interval and is incorporated into the controller before execution. Subsequently, \cite{wang2024Target} provides a more detailed analysis of the sampling-induced SSE and the parameter conditions arising from discretization, but does not develop an SSE cancellation mechanism.

Despite these advances, the discrete-time BoTLC problem remains underexplored. Therefore, the aforementioned discussion motivates the present paper to develop an adaptive alternative to static correction mechanism in \cite{chen2022DiscreteTime}. The main contributions are summarized as follows.
\begin{itemize}
    \item We propose a novel adaptive compensator for the discrete-time circumnavigation controller and establish exponential convergence. Driven by the estimated agent-to-target distance, the proposed compensator dynamically learns the geometric correction required to eliminate the sampling-induced orbit bias. In contrast to the static correction in \cite{chen2022DiscreteTime}, neither the controller nor the compensator update explicitly uses the sampling interval. The proposed approach also goes beyond the SSE characterization in \cite{wang2024Target} by providing a mechanism that eliminates this error.

    \item We provide an explicit characterization of the persistence of excitation (PE) of the bearing measurements, thereby directly revealing how the excitation strength is related to the tangential step and the distance bound. This, in turn, enables an explicit PE-based estimate of the exponential decay rate for the estimation error, whereas \cite{chen2022DiscreteTime} relies on an implicit existence claim and \cite{wang2024Target} uses geometric contraction without a quantitative construction using the PE constants.

\end{itemize}

\section{PRELIMINARIES}
\label{sec:Preliminaries}
        
    \subsection{Background and Notation}
    \label{sec:background_and_notations}
        
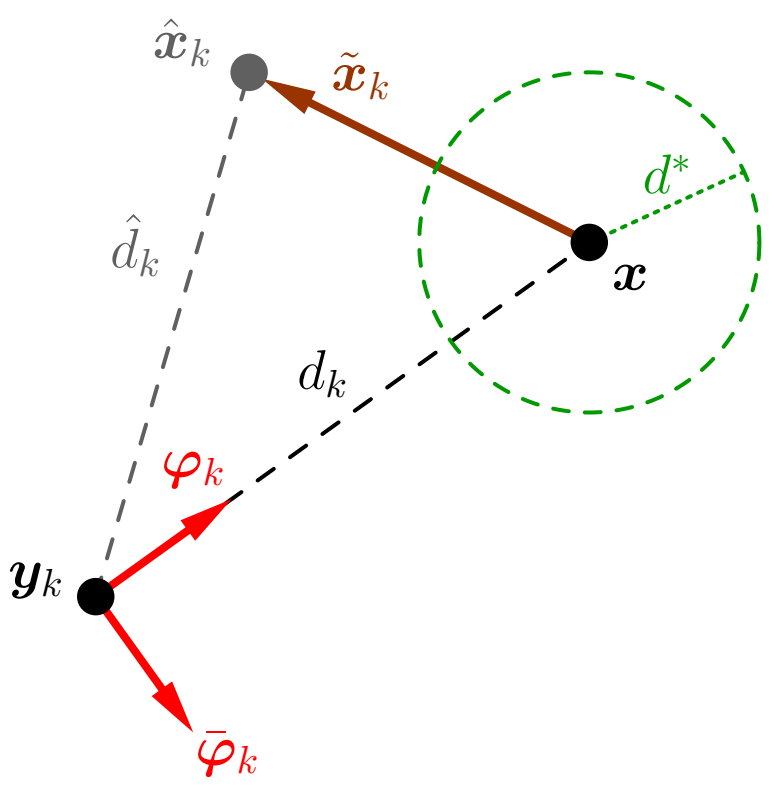
\begin{figure}[t]
            \centering
            \scalebox{0.75}{
                \begin{tikzpicture}[scale = 1.5, line cap=round,line join=round,>=triangle 45,x=1cm,y=1cm]
                    \clip(-4,-3.4) rectangle (1.2,1.56);
                    %
                    %
                    \draw [line width=1pt,dash pattern=on 5pt off 5pt] (-2.9019083532629173,-2.0826434806438683)-- (0,0);
                    \draw (-1.8,-0.55) node[anchor=north west] {\Large $d_{k}$};
                    %
                    %
                    \draw [line width=1pt,dash pattern=on 5pt off 5pt,color=wqwqwq] (-2.9019083532629173,-2.0826434806438683)-- (-2,1);
                    \draw [color=wqwqwq](-2.9,0.25) node[anchor=north west] {\Large $\hat{d}_{k}$};
                    
                    %
                    %
                    \draw [-{Triangle[length=5mm,width=2mm]},line width=2pt,color=ffqqqq] (-2.9019083532629173,-2.0826434806438683) -- (-2.0894814568450726,-1.499580415464064);
                    \draw [color=ffqqqq](-2.6,-1.15) node[anchor=north west] {\Large $\bm{\varphi}_{k}$};
                    %
                    %
                    \draw [-{Triangle[length=5mm,width=2mm]},line width=2pt,color=ffqqqq] (-2.9019083532629173,-2.0826434806438683) -- (-2.318845288083113,-2.895070377061713);
                    \draw [color=ffqqqq](-2.4,-2.8) node[anchor=north west] {\Large $\bar{\bm{\varphi}}_{k}$};
                    \draw [fill=black] (-2.9019083532629173,-2.0826434806438683) circle (3pt);
                    \draw (-3.5,-1.8) node[anchor=north west] {\Large $\bm{y}_{k}$};

                    \draw [-{Triangle[length=5mm,width=2mm]},line width=2pt,color=zzttqq] (0,0) -- (-1.94, 0.97);
                    \draw [color=zzttqq](-1.6,1.2) node[anchor=north west] {\Large $\tilde{\bm{x}}_{k}$};
        
                    \draw [line width=1pt,dash pattern=on 4pt off 4pt,color=qqzzqq] (0,0) circle (1cm);
                    \draw [line width=1pt,dotted,color=qqzzqq] (0,0)-- (0.9097994775058879,0.41504808243143776);
                    \draw [color=qqzzqq](0.23133379900014597,0.6) node[anchor=north west] {\Large $d^*$};
        
                    \draw [fill=black] (0,0) circle (3pt);
                    \draw (0.05,-0.05) node[anchor=north west] {\Large $\bm{x}$};

                    \draw [draw = wqwqwq, fill=wqwqwq] (-2,1) circle (3pt);
                    \draw [color=wqwqwq](-2.65,1.4) node[anchor=north west] {\Large $\hat{\bm{x}}_{k}$};
        


                \end{tikzpicture}
             }
            \caption{\label{fig:notation_conventions}Problem geometry and notations. Symbols $\hat{d}_k$, $\hat{\bm{x}}_k$, $\tilde{\bm{x}}_k$, and $d^*$ are included for completeness and will be explained in the following sections.}
        \end{figure}

        Consider a planar setting that contains a stationary target and an agent tasked with the BoTLC mission. Let $\bm{y}_k:=\bm{y}(kT)\in\mathbb{R}^2$ be the position of the agent at time $kT$ where $k\in\mathbb{Z}_{\geq 0}$ and $T>0$ is the sampling and control-update interval. Similarly, let $\bm{x}_{k}:=\bm{x}(kT)\in\mathbb{R}^2$ be the position of the target at time $kT$, and we suppose that this quantity is \emph{unknown a priori} to the agent. Since the target is assumed stationary, we drop the subscript and write $\bm{x}_k\equiv \bm{x}$ for all $k\geq 0$. We further suppose that the agent obeys the following first-order kinematics:
        \begin{equation}
            \label{eq:agent_kinematics_DT}
            \bm{y}_{k+1} = \bm{y}_{k}+\bm{u}_{k}T, \quad k\in\mathbb{Z}_{\geq 0},
        \end{equation}
        where $\bm{u}_k:=\bm{u}(kT)\in\mathbb{R}^2$ is the control signal to be designed. Consequently, the agent follows a rectilinear trajectory determined by the velocity $\bm{u}_k$ during the time interval from $kT$ to $(k+1)T$. The available bearing measurement is represented by a unit vector $\bm{\varphi}_k$, defined by
        \begin{equation}
            \label{eq:def_varphi}
            \bm{\varphi}_k =\frac{\bm{x}-\bm{y}_{k}}{\norm{\bm{x}-\bm{y}_{k}}} =: \frac{\bm{x}-\bm{y}_{k}}{d_{k}},
        \end{equation}
        where $\|\cdot\|$ denotes the Euclidean norm and $d_k$ represents the distance between the agent and the target at time $kT$. Let $\bar{\bm{\varphi}}_{k}\in\mathbb{R}^2$ be a unit vector perpendicular to $\bm{\varphi}_{k}$, obtained by $\pi/2$ clockwise rotation of $\bm{\varphi}_{k}$. See Fig.~\ref{fig:notation_conventions} for illustration.

    \subsection{Problem Statements}
    \label{sec:problem_statements}
        \begin{problem}[Target Localization]
            \label{problem:Localization_Problem}
            Let $\hat{\bm{x}}_k\in\mathbb{R}^2$ denote the agent's estimate of the target position. Design a discrete-time target estimator $\hat{\bm{x}}_{k+1}=\bm{f}\big(\hat{\bm{x}}_k, \bm{y}_k,\bm{\varphi}_k\big)$ to localize a stationary target $\bm{x}$. Specifically, the position estimation error $\tilde{\bm{x}}_k$, defined as,
            \begin{equation}
                \label{eq:def_x_tilde}
                \tilde{\bm{x}}_k:=\hat{\bm{x}}_k-\bm{x},
            \end{equation}
            must converge to zero exponentially fast as $k\to\infty$.
        \end{problem}

        \vspace{0.25cm}
        
        \begin{problem}[Target Circumnavigation]
            \label{problem:Circumnavigation_Problem}
            Let $d^*>0$ denote the radius of the desired circumnavigation orbit. Design a discrete-time circumnavigation controller in the form of $\bm{u}_{k}=\bm{g}\big(\bm{y}_k,\bm{\varphi}_k,\hat{\bm{x}}_k,d^*\big)$, without explicitly using the sampling interval $T$ in the control update, to ensure that the sampled circumnavigation tracking error, defined as
            \begin{equation}
                \label{eq:def_tracking_error}
                \delta_k := \norm{\bm{x}  - \bm{y}_k} - d^{*},
            \end{equation}
            converges to zero exponentially fast as $k\to\infty$.
        \end{problem}

    \subsection{Motivation: Shortcomings of Existing Methods}
    \label{sec:Shortcomings of Existing Methods}

        In this section, we demonstrate that a direct discretization of the estimator-coupled control algorithms initially proposed in \cite{deghat2010Target,deghat2014Localization} fails to solve \textit{Problem~\ref{problem:Circumnavigation_Problem}}. Specifically, we show that the directly discretized controller  introduces a steady-state tracking error $\delta_{\infty}>0$ such that $\delta_k\to\delta_\infty \neq 0$ as $k\to\infty$. To this end, consider the below discrete-time formulation of the continuous-time control and estimation algorithms in \cite{deghat2010Target,deghat2014Localization},
        \begin{align*}
            \bm{u}_k &= k_c \Big(\hat{d}_k-d^*\Big) \bm{\varphi}_k+k_\omega \bar{\bm{\varphi}}_k , \numberthis \label{eq:def_original_controller_DT} \\
            \hat{d}_k &= \|\hat{\bm{x}}_k-\bm{y}_k\|,  \numberthis \label{eq:def_d_hat_DT} \\
            \hat{\bm{x}}_{k+1} &= \hat{\bm{x}}_k+ k_x \Big(\bm{I}_2 - \bm{\varphi}_k\bm{\varphi}^{\top}_{k}\Big)\big(\bm{y}_k-\hat{\bm{x}}_k\big), \numberthis \label{eq:def_target_estimator_DT}
        \end{align*}
        where $k_c,k_\omega>0$ are control gains, $k_x\in(0,2)$ is the estimator gain\footnote{The selection range of the estimator gain $k_x$ will be explained later in the convergence analysis.}, and $\bm{I}_2$ denotes a $2$-by-$2$ identity matrix. 

        For notational brevity, let $\tau_c:=k_cT$ and $\tau_\omega:=k_\omega T$, and denote the distance estimation error by
        \begin{equation}
            \label{eq:def_varrho}
            \varrho_k:=d_k-\hat{d}_k.
        \end{equation}
        Then, the evolution of the agent-target distance $d_k$, defined in \eqref{eq:def_varphi}, is described by the below difference equation,
        \begin{align*}
            d_{k+1} \overset{\eqref{eq:def_varphi}}&{=} \norm{\bm{x}-\bm{y}_{k+1}} \overset{\eqref{eq:agent_kinematics_DT}}{=} \norm{\bm{x}-(\bm{y}_k+\bm{u}_k \cdot T)} \\
            \overset{\eqref{eq:def_original_controller_DT}}&{=} \left\| d_k \bm{\varphi}_k-\bigl(k_c(\hat{d}_k-d^*)\bm{\varphi}_k+k_\omega \bar{\bm{\varphi}}_k\bigr)T \right\| \\
            \overset{\eqref{eq:def_varrho}}&{=} \left\|
            \begin{bmatrix}
                (1-\tau_c)d_k+\tau_c(\varrho_k+d^*)\\
                \tau_\omega
            \end{bmatrix}
            \right\|  =:h_1(d_k,\varrho_k), \numberthis \label{eq:def_d_k+1_original}
        \end{align*}
        where in the last line of \eqref{eq:def_d_k+1_original}, we used the fact that, for perpendicular unit vectors $\bm{\varphi}_k$ and $\bar{\bm{\varphi}}_k$, it holds for any arbitrary $\alpha$ and $\beta$ that $\|\alpha \bm{\varphi}_k + \beta \bar{\bm{\varphi}}_k \| = \sqrt{\alpha^2 + \beta^2} = \|[\alpha, \; \beta]^\top\|$. We now proceed to show that even with the ideal initialization $\hat{\bm{x}}_0=\bm{x}$, under which the estimator \eqref{eq:def_target_estimator_DT} gives $\hat{\bm{x}}_k\equiv \bm{x}$ and hence $\varrho_k\equiv 0$ for all $k\in\mathbb{Z}_{\geq 0}$, the agent driven by the naively discretized controller \eqref{eq:def_original_controller_DT} is unable to achieve zero steady-state error in the circumnavigation tracking error.

        \begin{proposition}
        \label{prop:biased_circumnavigation_original}
        Suppose that $d^*>0$, $\tau_\omega>0$, and $0<\tau_c<2$. Then, the following autonomous iteration
        \begin{equation}
            \label{eq:def_h_1(d_k,0)}
            d_{k+1}=h_1(d_k,0), \quad k\in\mathbb{Z}_{\geq 0},
        \end{equation}
        has a unique positive \emph{fixed point} (in the sense of \cite[Def.~1.6]{elaydi2005Introduction})
        \begin{equation}
            \label{eq:def_d_infty_original}
            d_\infty =
            \frac{
                (1-\tau_c)d^*
                +\sqrt{(d^*)^2+
                    \frac{2-\tau_c}{\tau_c}\tau_\omega^2}
            }{2-\tau_c}
        \end{equation}
        that is globally exponentially stable (GES) on $\mathbb{R}_{>0}$. Specifically, for every $d_0>0$, the agent-target distance satisfies
        \begin{equation}
            \label{eq:d_convergence_original}
            |d_k-d_\infty|
            \leq |1-\tau_c|^k\, |d_0-d_\infty|,
            \quad \forall\, k\geq 1.
        \end{equation}
        Consequently, the tracking error $\delta_k=d_k-d^*$ converges exponentially to $\delta_\infty:=d_\infty-d^*>0$.
    \end{proposition}

    \begin{proof}
        Under $\varrho_k\equiv 0$, the map \eqref{eq:def_d_k+1_original} gives, for any $d>0$,
        \begin{equation}
            \label{eq:original_dis_map_zero_est_error}
            h_1(d,0)=
            \left\|
            \begin{bmatrix}
                (1-\tau_c)d+\tau_c d^*\\
                \tau_\omega
            \end{bmatrix}
            \right\|.
        \end{equation}
        Since $h_1(d,0)\geq\tau_\omega>0$, the iteration preserves positivity. Solving $h_1(d,0)=d$ for $d>0$ gives the fixed point in \eqref{eq:def_d_infty_original}. For any $d,s>0$, applying the triangle inequality gives
        \begin{align*}
            \scalemath{0.95}{|h_1(d,0)-h_1(s,0)|
            \overset{\eqref{eq:original_dis_map_zero_est_error}}{\leq}
            \left\|
            \begin{bmatrix}
                (1-\tau_c)(d-s)\\
                0
            \end{bmatrix}
            \right\| =|1-\tau_c|\,|d-s|.}
        \end{align*}
        Because $\tau_c$ must be chosen to meet $\tau_c\in(0,2)$, we have $|1-\tau_c|<1$. Hence the map \eqref{eq:def_h_1(d_k,0)} is contractive and therefore cannot admit two distinct positive fixed points. Taking $s=d_\infty$, we obtain
        \begin{equation*}
            |d_{k+1}-d_\infty|
            \leq |1-\tau_c|\,|d_k-d_\infty|.
        \end{equation*}
        Iterating this inequality establishes \eqref{eq:d_convergence_original} and hence GES. In particular, when $\tau_c=1$, the fixed point is reached in one step. Now, subtracting $d^*$ from \eqref{eq:def_d_infty_original} gives
        \begin{equation*}
            \delta_\infty=
            \frac{
                \sqrt{(d^*)^2+
                    \frac{2-\tau_c}{\tau_c}\tau_\omega^2}
                -d^*
            }{2-\tau_c}>0,
        \end{equation*}
        where the strict inequality follows from the assumption that $0<\tau_c<2$ and $\tau_\omega>0$. By construction, $|\delta_k-\delta_\infty|=|d_k-d_\infty|$, and hence the tracking error $\delta_k$ converges exponentially to this strictly positive limit.
    \end{proof}

        \begin{remark}
            The condition $0<\tau_c<2$ in \emph{Prop.~\ref{prop:biased_circumnavigation_original}} ensures that $|1-\tau_c|<1$ and hence provides a sufficient condition for the GES to the unique positive fixed point $d_\infty$. This condition is not, however, necessary for the existence of such a positive fixed point. For fixed $k_c,k_\omega,d^*>0$, substituting $\tau_c=k_cT$ and $\tau_\omega=k_\omega T$ into \eqref{eq:def_d_infty_original} and expanding allow us to obtain
            \begin{equation}
                \label{eq:small_T_bias_original}
                \delta_\infty
                = \frac{k_\omega^2}{2k_c d^*}\,T+O(T^2),
                \quad T\to 0.
            \end{equation}
            Thus, the steady-state tracking error is first order in the sampling period for sufficiently small $T$, and $d_\infty\to d^*$ as $T\to 0$ implies that no SSE in radial tracking exists when the system is in continuous time. This implication is also consistent with classical continuous-time analyses in \cite[Lem.~1]{deghat2010Target} and \cite[Lem.~1]{deghat2014Localization}.
        \end{remark}

\begin{figure}[t]
            \centering
            \begingroup%
\definecolor{seqBlue}{RGB}{0,0,255}%
\definecolor{seqRed}{RGB}{255,0,0}%
\definecolor{seqOrange}{RGB}{255,127,0}%
\definecolor{seqGreen}{RGB}{0,153,0}%
\definecolor{seqMagenta}{RGB}{255,0,255}%
\definecolor{seqPurple}{RGB}{127,0,255}%
\pgfmathsetlengthmacro{\seqUnit}{\linewidth/18}%
\begin{tikzpicture}[
  x=\seqUnit,y=\seqUnit,
  line cap=round,line join=round,
  font=\fontsize{8}{10}\selectfont,
  every node/.style={inner sep=0pt,outer sep=0pt},
  seq line/.style={line width=.7pt},
  seq guide/.style={line width=.4pt,dash pattern=on 2pt off 1.6pt},
  seq vector/.style={seqRed,line width=1pt,
    -{Triangle[length=4pt,width=3pt]}},
  seq open arrow/.style={-{Straight Barb[length=2.8pt,width=3.5pt]}}
]

\path (-11,-6) rectangle (4.75,3);

\coordinate (X) at (0,0);
\coordinate (Y) at (-7.071068,0);
\coordinate (Ynext) at (-7.071068,-5.824376);
\coordinate (Yprime) at (-4.009570,-5.824376);
\coordinate (Foot) at (Yprime |- X);
\coordinate (Radial) at (-5.457938,-4.495655);
\coordinate (DimLeft) at ($(Y)+(0,.8)$);
\coordinate (DimRight) at (Foot |- DimLeft);

\begin{scope}
  \clip (-7.3,-7.2) rectangle (5.1,2.95);
  \draw[line width=.8pt] (-45:7.071068)
    arc[start angle=-45,end angle=-225,radius=7.071068];
\end{scope}

\draw[seq guide] (Y) -- (X);
\draw[seq guide] (Ynext) -- (X);
\draw[seq guide] (Y) -- (DimLeft);
\draw[seq guide] (Foot) -- (DimRight);
\draw[seq line,seqOrange] (Yprime) -- (X);
\draw[seq line,seqGreen] (Y) -- (Ynext);
\draw[seq line,seqGreen,seq open arrow] (Yprime) -- (Foot);
\draw[seq line,seqMagenta] (Ynext) -- (Yprime);
\draw[seq line,seqMagenta,<->,
  >={Straight Barb[length=2.8pt,width=3.5pt]}] (DimLeft) -- (DimRight);
\draw[seqPurple,line width=.9pt,seq open arrow] (Ynext) -- (Radial);

\draw[seq vector] (Y) -- ++(2,0);
\draw[seq vector] (Y) -- ++(0,-2);

\node[seqBlue,anchor=east] at ([xshift=-5pt]Ynext)
  {$\boldsymbol{y}_{k+1}$};
\node[seqBlue,anchor=south west] at ([xshift=5pt,yshift=-3pt]Yprime)
  {$\boldsymbol{y}_{k+1}'$};
\node[seqBlue,anchor=south east] at ([xshift=-4pt,yshift=5pt]Y)
  {$\boldsymbol{y}_k$};
\node[anchor=south west] at ([xshift=4pt,yshift=4pt]X)
  {$\boldsymbol{x}$};
\node[seqRed,anchor=north] at (-5.40,-.28)
  {$\boldsymbol{\varphi}_k$};
\node[seqRed,anchor=east] at (-7.38,-1.65)
  {$\bar{\boldsymbol{\varphi}}_k$};
\node[seqOrange,anchor=west] at (-1.80,-3.05) {$d^*$};
\node[seqPurple,anchor=south] at (-6.35,-4.7) {$\delta_{k+1}$};

\node[seqGreen,anchor=east] at (-7.45,-3.20) {$
  \begin{aligned}
    \tau_\omega &=\lVert k_\omega\bar{\boldsymbol{\varphi}}_k\rVert\cdot T\\
    &=k_\omega T
  \end{aligned}$};
\node[seqMagenta,anchor=south west] at (-6.25,1.30) {$
  \begin{aligned}
    \tau_c \zeta_{\infty} &= d^* - \sqrt{(d^*)^2-\tau_\omega^2}
  \end{aligned}$};

\filldraw[draw=black,fill=seqBlue,line width=.35pt]
  (Y) circle[radius=2pt];
\filldraw[draw=black,fill=seqBlue,line width=.35pt]
  (Yprime) circle[radius=2pt];
\draw[draw=seqBlue,fill=white,line width=1.1pt]
  (Ynext) circle[radius=2.1pt];
\fill[black] (X) circle[radius=2pt];
\end{tikzpicture}%
\endgroup%
            \caption{\label{fig:geometric_bias}One-step geometry of the sampling-induced circumnavigation tracking error and the required inward correction.}
        \end{figure}
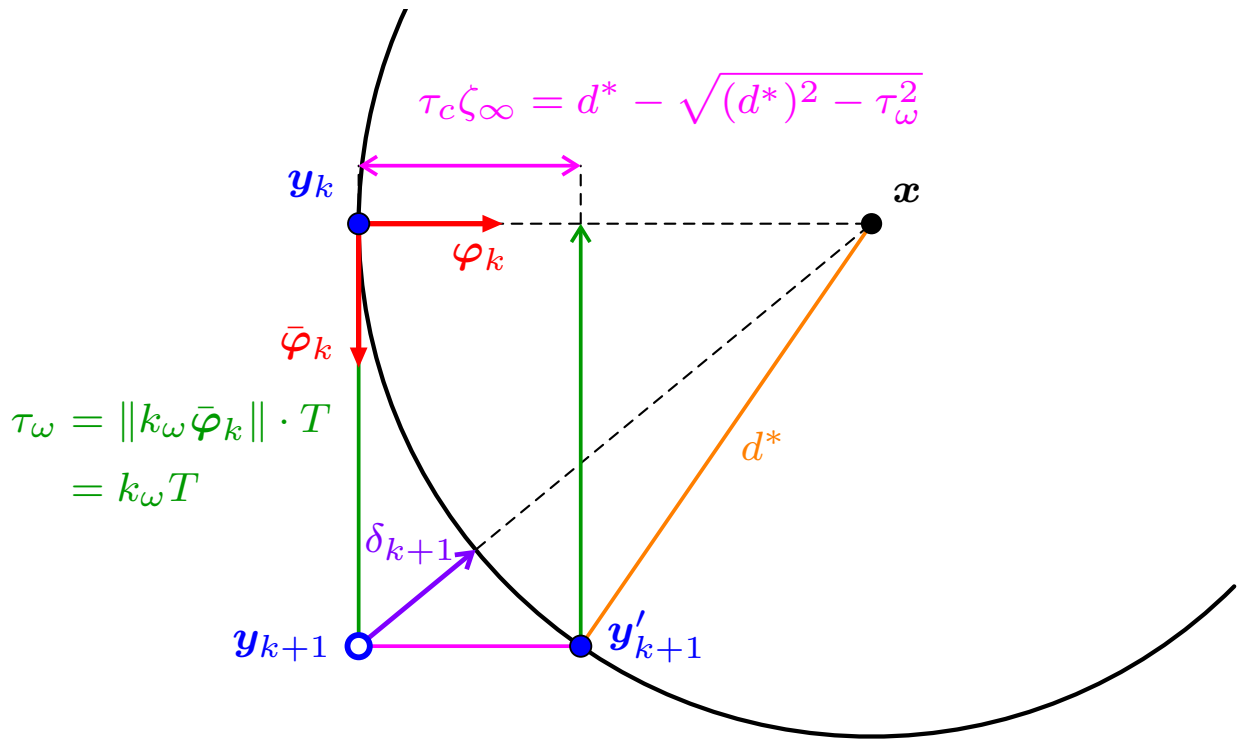

        \begin{remark}
            \label{rem:geometric_bias_original}
            The steady-state error discussed in \emph{Prop~\ref{prop:biased_circumnavigation_original}} also admits a geometric interpretation. Figure~\ref{fig:geometric_bias} illustrates the origin of the tracking bias. Under exact localization $\hat{\bm{x}}_k=\bm{x}$ and when $d_k=d^*$, the naive controller \eqref{eq:def_original_controller_DT} produces a purely tangential displacement of length $\tau_\omega$, moving the agent from $\bm{y}_k$ to $\bm{y}_{k+1}$ rather than the intended $\bm{y}'_{k+1}$. This leads to
            \begin{equation}
                \label{eq:tangential_step_bias_original}
                d_{k+1}=\sqrt{(d^*)^2+\tau_\omega^2}>d^*.
            \end{equation}
            For $0<\tau_\omega<d^*$, reaching the depicted point $\bm{y}'_{k+1}$ on the circle with the same tangential displacement requires the additional inward displacement $\tau_c\zeta_{\infty}$ shown in Fig.~\ref{fig:geometric_bias}, where, for $0<k_\omega T<d^*$, direct trigonometric calculation gives
            \begin{equation}
                \label{eq:def_zeta_ss_geometric}
                \zeta_{\infty}
                :=
                \frac{
                    d^*-\sqrt{(d^*)^2-(k_\omega T)^2}
                }{k_cT}.
            \end{equation}
            Therefore, it is also evident from a geometric argument that a nonzero inward correction is needed even when the sampled tracking error vanishes.
        \end{remark}

\section{PROPOSED ALGORITHM}
\label{sec:proposed_algorithm}

    \begin{subequations}
        \label{eq:def_controller_integral_DT}
        To eliminate the SSE in tracking, we propose the following improved circumnavigation control law with integral action,
        \begin{equation}
            \label{eq:def_circum_controller_DT}
                \bm{u}_k = k_c \Big( \hat{d}_k -d^*+\hat{\zeta}_k\Big)\bm{\varphi}_k+k_\omega \bar{\bm{\varphi}}_k,
        \end{equation}
        where $\hat{d}_k$ is defined in \eqref{eq:def_d_hat_DT} with $\hat{\bm{x}}_k$ updated by \eqref{eq:def_target_estimator_DT}, and $\hat{\zeta}_k$ is the proposed adaptive compensator, 
        \begin{align}
            \label{eq:def_integral_term_DT}
            \hat{\zeta}_{k+1} = d^* - \frac{d^*-\hat{\zeta}_k}{1+k_\zeta \bigl(\hat{d}_k-d^*\bigr)/d^*}, \quad \hat{\zeta}_0 := 0,
        \end{align}
        and $k_\zeta \in (0,1)$ is a compensator gain. Note that neither the controller \eqref{eq:def_circum_controller_DT} nor the compensator update \eqref{eq:def_integral_term_DT} requires explicit knowledge of the sampling interval $T$.
    \end{subequations}

    \begin{remark}
        The proposed compensator \eqref{eq:def_integral_term_DT} admits an interpretation in terms of an adaptive reference radius. Let
        \begin{equation}
            \label{eq:def_reference_radius_DT}
            r_k:=d^*-\hat{\zeta}_k.
        \end{equation}
        Then, the controller \eqref{eq:def_circum_controller_DT} can be written as $\bm{u}_k=k_c(\hat{d}_k-r_k)\bm{\varphi}_k +k_\omega\bar{\bm{\varphi}}_k$, and it follows from \eqref{eq:def_integral_term_DT} that
        \begin{equation}
            \label{eq:reference_radius_update_DT}
            r_{k+1}
            =\frac{r_k}{1+k_\zeta(\hat{d}_k-d^*)/d^*},
            \quad r_0=d^*.
        \end{equation}
        By \eqref{eq:def_d_hat_DT}, $\hat{d}_k\geq0$, and since $0<k_\zeta<1$, the denominator in \eqref{eq:reference_radius_update_DT} is bounded below by $1-k_\zeta>0$. Consequently, $r_k>0$, or equivalently $\hat{\zeta}_k<d^*$, for all $k\in\mathbb{Z}_{\geq0}$. Also note that the adaptive reference radius $r_k$ decreases when $\hat{d}_k>d^*$ and increases when $\hat{d}_k<d^*$. Further, by subtracting $\hat{\zeta}_k$ from both sides of \eqref{eq:def_integral_term_DT},  the compensator update can be expressed as the following difference equation,
        \begin{equation}
            \label{eq:nonlinear_integral_increment_DT}
            \hat{\zeta}_{k+1}-\hat{\zeta}_k
            =\frac{k_\zeta(d^*-\hat{\zeta}_k)}
            {(1-k_\zeta)d^*+k_\zeta\hat{d}_k}
            \Big(\hat{d}_k-d^*\Big).
        \end{equation}
        Note that the coefficient multiplying the estimated tracking error $\hat{d}_k-d^*$ in \eqref{eq:nonlinear_integral_increment_DT} is always positive. Thus, the compensator increment $\hat{\zeta}_{k+1}-\hat{\zeta}_k$ vanishes \emph{if and only if} $\hat{d}_k=d^*$.
    \end{remark}

\section{STABILITY ANALYSIS}
\label{sec:stability_analysis}
    
    Before we proceed, we impose the following assumptions.

    \begin{assumption}
        \label{assumption:d(0)_initial_condition}
        We assume that the agent and the target occupy different locations at $t=0$, that is, $d_0>0$. 
    \end{assumption}

    \begin{assumption}
        \label{assumption:sampling_time}
        Let $d^*,k_c,k_\omega,T>0$, $k_x\in(0,2)$, $k_\zeta\in(0,1)$, and recall that $\tau_c=k_c T$ and $\tau_\omega = k_\omega T$. We suppose that the sampling interval and gains satisfy
        \begin{equation}
            \label{eq:radial_feasibility_DT}
            0<\tau_c\leq1,\quad
            0<\tau_\omega<d^*\sqrt{\tau_c(2-\tau_c)}.
        \end{equation}
    \end{assumption}

    Under the standing \textit{Asms.~\ref{assumption:d(0)_initial_condition}-\ref{assumption:sampling_time}}, substituting \eqref{eq:def_circum_controller_DT} into \eqref{eq:agent_kinematics_DT} gives, for all $k \in\mathbb{Z}_{\geq 0}$,
    \begin{equation}
        \label{eq:def_d_k+1_integral} 
        d_{k+1} = \sqrt{\bigl((1-\tau_c)d_k + \tau_c(\varrho_k+r_k)\bigr)^2+\tau_{\omega}^{2}},
    \end{equation}
    which implies $d_{k+1}\geq\tau_\omega>0$.

    \begin{lemma}
        \label{lemma:x_tilde_is_non_increasing_DT}
        Suppose \textit{Asms.~\ref{assumption:d(0)_initial_condition}-\ref{assumption:sampling_time}} hold. Then, under the target estimator \eqref{eq:def_target_estimator_DT} and the controller with integral action \eqref{eq:def_controller_integral_DT}, the position estimation error $\tilde{\bm{x}}_k$, defined in \eqref{eq:def_x_tilde}, satisfies the below inequality bound:
        \begin{equation}
            \label{eq:x_tilde_is_non_increasing}
            |\varrho_k| \leq \norm{\tilde{\bm{x}}_k}\leq \norm{\tilde{\bm{x}}_0}, \quad \forall \, k \in\mathbb{Z}_{\geq 0},
        \end{equation}
        where $\varrho_k$ is the distance estimation error defined in \eqref{eq:def_varrho}. 
    \end{lemma}

     \begin{proof}
        From \cite[Lem.~2]{deghat2010Target}, we obtain the following relationship in discrete-time, 
        \begin{equation}
            \label{eq:projection_property}
            \bm{I}_2-\bm{\varphi}_k\bm{\varphi}_{k}^{\top} = \bar{\bm{\varphi}}_k \bar{\bm{\varphi}}_k^\top := P_k.
        \end{equation}
        Observe that $P_{k}^\top = P_{k}^{2}=P_{k}$, and bearing geometry gives $P_k(\bm{y}_k-\bm{x})=0$. The error dynamics of the target estimation error $\tilde{\bm{x}}_k$ therefore satisfies
        \begin{align}
            \label{eq:x_tilde_k+1}
            \tilde{\bm{x}}_{k+1}  = \bigl( I_2 - k_x P_k\bigr) \tilde{\bm{x}}_k .
        \end{align}
        Now consider the following Lyapunov candidate 
        \begin{equation}
            \label{eq:estimation_error_Lyapunov}
            V_k=\tilde{\bm{x}}_{k}^\top\tilde{\bm{x}}_{k}=\|\tilde{\bm{x}}_k\|^2,
        \end{equation}
        whose one-step difference along \eqref{eq:x_tilde_k+1} is, 
        \begin{align*}
            V_{k+1}-V_k &= \tilde{\bm{x}}_{k}^{\top} \bigl( (I_2 - k_x P_k)^2 - I_2 \bigr) \tilde{\bm{x}}_k \\
            &= - k_x (2- k_x) \bigl| \bar{\bm{\varphi}}_{k}^{\top}\tilde{\bm{x}}_k|^2 \leq 0, \numberthis \label{eq:est_Lyapunov_inequality}
        \end{align*}
        which is negative semi-definite given that \textit{Asm.~\ref{assumption:sampling_time}} holds. Thus, we obtain that $\tilde{\bm{x}}_k$ is bounded and satisfies $\|\tilde{\bm{x}}_k\| \leq \|\tilde{\bm{x}}_0\|$ for all $k\in\mathbb{Z}_{\geq 0}$. Using the triangle inequality, we further obtain that $|d_k-\hat{d}_k|=|\varrho_k| \leq \|\tilde{\bm x}_k\|$. The inequality in \eqref{eq:x_tilde_is_non_increasing} can be immediately obtained by combining the above two inequalities here.
    \end{proof}

    Since the signal $\bm{\varphi}_k$ is explicitly used in the target estimator \eqref{eq:def_target_estimator_DT} and the controller with integral action \eqref{eq:def_controller_integral_DT}, we need to ensure that $\bm{\varphi}_k$ is well-defined for all $k\in\mathbb{Z}_{\geq 0}$. We now establish bounds on the agent-target distance directly from the distance recursion and the compensator update.

\begin{lemma}
    \label{lemma:d_k_is_bounded_DT}
    Suppose \textit{Asms.~\ref{assumption:d(0)_initial_condition} and~\ref{assumption:sampling_time}} hold. Then, under the target estimator \eqref{eq:def_target_estimator_DT} and the controller with integral action \eqref{eq:def_controller_integral_DT}, the adaptive reference radius $r_k$, the agent-target distance $d_k$, and the compensator satisfy the following bounds,
    \begin{subequations}
        \begin{gather}
            0<r_k\leq r_{\max}, \label{eq:r_k_bound} \\
            0<d_{\min}:=\min\{d_0,\, \tau_\omega\} \leq d_k\leq d_{\max}, \label{eq:d_k_is_bounded_DT} \\
            d^*-r_{\max}\leq\hat{\zeta}_k<d^*, \label{eq:compensator_bound}
        \end{gather}
    \end{subequations}
    for all $k\in\mathbb{Z}_{\geq 0}$, where 
    \begin{subequations}
        \label{eq:direct_comparison_constants_DT}
        \begin{align}
            r_{\max}&:=\frac{
                \max\left\{
                    (1+\tau_c)\|\tilde{\bm{x}}_0\|,\;
                    \frac{(1-k_\zeta)d^*}{k_\zeta}
                \right\}
            }{\tau_c(1-k_\zeta)^2}, \label{eq:def_r_max} \\
            d_{\max}&:=\max\left\{
                d_0, \;
               \|\tilde{\bm{x}}_0\|+r_{\max} + \frac{\tau_\omega}{\tau_c}
                \right\}. \label{eq:def_d_max}
        \end{align}
    \end{subequations}
\end{lemma}

\vspace{0.1cm}

\begin{proof}
    Since the distance recursion \eqref{eq:def_d_k+1_integral} implies that $d_{k+1}\geq\tau_\omega>0$, $d_0>0$ by \emph{Asm.~\ref{assumption:d(0)_initial_condition}}, and the denominator in the $r_k$ update \eqref{eq:reference_radius_update_DT} is at least $1-k_\zeta >0$, the sampled bearing $\bm{\varphi}_k$ and hence the controller updates remain well-defined. The lower bound in \eqref{eq:d_k_is_bounded_DT} then follows immediately.

    Next, we establish the claimed upper bound on $r_k$ in \eqref{eq:r_k_bound}. From \textit{Lem.~\ref{lemma:x_tilde_is_non_increasing_DT}}, we know that $|\varrho_k|\leq \|\tilde{\bm{x}}_0\|$ for all $k\in\mathbb{Z}_{\geq0}$. Applying \eqref{eq:def_varrho} and using the distance recursion \eqref{eq:def_d_k+1_integral}, together with $0<\tau_c\leq 1$, gives
    \begin{equation*}
        \hat{d}_{k+1} \geq d_{k+1} - \|\tilde{\bm{x}}_0\| \geq \tau_c r_k - E, \quad E:=(1+\tau_c)\|\tilde{\bm{x}}_0\|.
    \end{equation*}
    For brevity, let $\xi:=1-k_\zeta$ and $\psi:=k_\zeta/d^*$. Since $\hat{d}_k,\hat{d}_{k+1}\geq0$ and $\tau_c r_k\leq\hat{d}_{k+1}+E$, applying \eqref{eq:reference_radius_update_DT} twice gives
    \begin{align*}
        \scalemath{0.975}{r_{k+2}}
        &= \scalemath{0.975}{\frac{r_k}
        {(\xi+\psi\hat{d}_k)(\xi+\psi\hat{d}_{k+1})}
        \leq \frac{\hat{d}_{k+1}+E}
        {\tau_c \xi(\xi+\psi\hat{d}_{k+1})}} \\
        &\leq \scalemath{0.975}{\frac{\max\{E/\xi,\,1/\psi\}}{\xi \tau_c}
        \overset{\eqref{eq:def_r_max}}{=}r_{\max}, }
    \end{align*}
    where the last inequality uses the fact that, for any $x\geq 0$, $x+E\leq\max\{E/\xi,1/\psi\}(\xi+\psi x)$. Since $r_0=d^*$ and $r_1\leq d^*/(1-k_\zeta)\leq r_{\max}$, \eqref{eq:r_k_bound} holds for all $k\in\mathbb{Z}_{\geq0}$.

    Finally, using \eqref{eq:def_d_k+1_integral} and the bounds on $\varrho_k$ and $r_k$, we obtain
    \begin{align*}
        d_{k+1}\leq(1-\tau_c)d_k +\tau_c(\|\tilde{\bm{x}}_0\|+r_{\max})+\tau_\omega.
    \end{align*}
    Then, standard scalar comparison gives $d_k\leq d_{\max}$. The remaining bounds follow from $d_{k+1}\geq\tau_\omega$ and $\hat\zeta_k=d^*-r_k$.
\end{proof}

    The distance upper bound $d_k\leq d_{\max}$ established in \eqref{eq:d_k_is_bounded_DT} of \emph{Lem.~\ref{lemma:d_k_is_bounded_DT}}, together with the nonzero tangential displacement $\tau_\omega$, provides a sufficient condition for persistence of excitation. In particular, it is sufficient to consider two consecutive bearing measurements. We formally establish this claim in the following lemma.

\begin{lemma}
    \label{lemma:PE_DT}
    Suppose \emph{Asms.~\ref{assumption:d(0)_initial_condition} and~\ref{assumption:sampling_time}} hold. Then, under the target estimator \eqref{eq:def_target_estimator_DT} and the controller with integral action \eqref{eq:def_controller_integral_DT}, the signal $\bar{\bm{\varphi}}_k$ is \emph{persistently exciting (PE)}. Specifically,
    \begin{equation}
        \label{eq:PE_inequalities}
        \sum_{j=k}^{k+1}
        \bar{\bm{\varphi}}_j\bar{\bm{\varphi}}_j^\top
        \geq \mu\bm{I}_2, \quad \mu:=1-\sqrt{1-
        \left(\frac{\tau_\omega}{d_{\max}}\right)^2}>0,
    \end{equation}
    holds for all $k\in\mathbb{Z}_{\geq 0}$ with $d_{\max}$ defined in \eqref{eq:def_d_max}. Consequently, the estimation error $\tilde{\bm{x}}_k$ converges to the origin exponentially fast.
\end{lemma}

\begin{proof}
    Since the target is stationary, the relative position update follows from \eqref{eq:agent_kinematics_DT}, \eqref{eq:def_varphi}, and \eqref{eq:def_circum_controller_DT} as
    \begin{equation}
        \label{eq:relative_position_update_DT}
            d_{k+1}\bm{\varphi}_{k+1}
            =\bigl((1-\tau_c)d_k
            +\tau_c(\varrho_k+r_k)\bigr)\bm{\varphi}_k -\tau_\omega\bar{\bm{\varphi}}_k.
    \end{equation}
    Premultiplying \eqref{eq:relative_position_update_DT} by $\bar{\bm{\varphi}}_k^\top$ and using $\bar{\bm{\varphi}}_k^\top\bm{\varphi}_k\equiv 0$ gives
    \begin{equation}
        \label{eq:consecutive_bearings_DT}
        \bar{\bm{\varphi}}_k^\top\bm{\varphi}_{k+1}
        =-\frac{\tau_\omega}{d_{k+1}}.
    \end{equation}
    Before proceeding, recall that for any pair of unit vectors $\bm{a},\bm{b}\in\mathbb{R}^2$, the eigenvalues of $\bm{a}\bm{a}^\top+\bm{b}\bm{b}^\top$ are $1\pm |\bm{a}^\top\bm{b}|$. Since the vectors $\{\bm{\varphi}_k,\bar{\bm{\varphi}}_k\}$ form an orthonormal basis of $\mathbb{R}^2$, we obtain
    \begin{align*}
        \lambda_{\min}\left(
            P_k + P_{k+1}    
        \right)
        &= 1- \sqrt{1-\left(\bar{\bm{\varphi}}_{k}^{\top}\bm{\varphi}_{k+1}\right)^2} \\
        \overset{\eqref{eq:consecutive_bearings_DT}}&{=} 1 - \sqrt{ 1 - \left( \frac{\tau_\omega}{d_{k+1}} \right)^2} \overset{\eqref{eq:d_k_is_bounded_DT}}{\geq} \mu > 0,
    \end{align*}
    where $P_j=\bar{\bm{\varphi}}_j\bar{\bm{\varphi}}_j^\top$ is defined earlier in \eqref{eq:projection_property}. Therefore, the signal $\bar{\bm{\varphi}}_k$ is PE.
    
    We now proceed to show that the estimation error $\tilde{\bm{x}}_k$ converges to the origin exponentially fast. Consider the Lyapunov candidate defined in \eqref{eq:estimation_error_Lyapunov} and let $p:=\bar{\bm{\varphi}}_k^\top\tilde{\bm{x}}_k$, $q:=\bar{\bm{\varphi}}_{k+1}^\top\tilde{\bm{x}}_{k+1}$, and $c:=\bar{\bm{\varphi}}_{k+1}^\top\bar{\bm{\varphi}}_k$. Then we obtain
    \begin{align*}
        \mu V_k \overset{\eqref{eq:x_tilde_k+1},\eqref{eq:est_Lyapunov_inequality},\eqref{eq:PE_inequalities}}&{\leq} p^2+(q+k_xcp)^2 \overset{|c|\leq 1}{\leq} 2(1+k_x^2)(p^2+q^2)\\
        \overset{\eqref{eq:est_Lyapunov_inequality}}&{=} \frac{2(1+k_x^2)}{k_x(2-k_x)} (V_k-V_{k+2}).
    \end{align*}
    It then follows immediately that
    \begin{equation}
        \label{eq:localization_two_step_contraction_DT}
        V_{k+2}\leq(1-\eta)V_k,\quad
        \eta:=\frac{k_x(2-k_x)\mu}{2(1+k_x^2)}\in(0,1),
    \end{equation}
    which completes the proof.
\end{proof}

\begin{theorem}
    \label{theorem:circum_error_goes_to_zero}
    Suppose \emph{Asms.~\ref{assumption:d(0)_initial_condition} and~\ref{assumption:sampling_time}} hold, and choose
    \begin{equation}
        \label{eq:radial_compensator_gain_DT}
        0<k_\zeta<\frac{\tau_c}{2+\tau_c}.
    \end{equation}
    Under the estimator \eqref{eq:def_target_estimator_DT} and the controller \eqref{eq:def_controller_integral_DT}, for every $d_0>0$ and $\hat{\bm{x}}_0\in\mathbb{R}^2$, it holds that
    \begin{equation}
        \label{eq:radial_convergence_integral_DT}
        d_k\to d^*,\quad
        r_k\to r_\infty:=
        \frac{\sqrt{(d^*)^2-\tau_\omega^2}-(1-\tau_c)d^*}
             {\tau_c}>0,
    \end{equation}
    as $k\to\infty$. Consequently, $\hat{\zeta}_k\to\zeta_\infty=d^*-r_\infty$, with $\zeta_\infty$ defined by \eqref{eq:def_zeta_ss_geometric}. Further, along every trajectory there exist $C>0$ and $\lambda\in(0,1)$, possibly depending on the initial condition, such that
    \begin{equation}
        \label{eq:radial_exponential_bound_DT}
        |d_k-d^*|+|r_k-r_\infty|\leq C\lambda^k,
        \quad \forall \, k\in \mathbb{Z}_{\geq 0}.
    \end{equation}
\end{theorem}

\begin{proof}
    This proof is organized into two parts: 
    \begin{itemize}
        \item \textbf{Part I} shows later in \eqref{eq:adaptive_reference_radius_bound_for_large_K} that there exists some sufficiently large but finite $K\geq 0$ such that the adaptive reference radius satisfies the eventual bound $r_k\leq 3d^*/2$ for all subsequent $k\geq K$. 

        \item \textbf{Part II} proceeds in three steps. \emph{Step 1} first introduces a coordinate transformation \eqref{eq:coordinate_transform} for the distance $d_k$ and reference radius $r_k$ updates. Then, by using the exact state-dependent coefficient (SDC) representation, these updates can be viewed as a forced system of the form $\bm{w}_{k+1}=M_k\bm{w}_k+\bm{\rho}_k$, defined later in \eqref{eq:radial_error_matrix_DT}, where $\bm{\rho}_k$ accounts for the forcing input due to the estimation error $\varrho_k$ and $\bm{w}_k$ is the transformed state. \emph{Step 2} and \emph{Step 3} use the polyhedral Lyapunov candidate \eqref{eq:polyhedral_Lyapunov} to establish the claimed exponential convergence bound \eqref{eq:radial_exponential_bound_DT}.
    \end{itemize}

    We now present the details of each part.

    \textbf{Part I.} First define the auxiliary variable
    \begin{equation}
        \label{eq:def_aux_s_k}
        s_k:= \gamma r_k -d_k, \quad \gamma := \frac{\tau_c (1-k_\zeta)}{(\tau_c(1-k_\zeta)+k_\zeta)}.
    \end{equation}
    The gain selection condition \eqref{eq:radial_compensator_gain_DT} implies that $\gamma>2/3$. The distance update \eqref{eq:def_d_k+1_integral} and compensator update \eqref{eq:def_integral_term_DT} imply
    \begin{align}
        \label{eq:d_and_r_inequalities}
        \begin{split}
            d_{k+1}&\geq(1-\tau_c)d_k+ \tau_c r_k-\tau_c |\varrho_k|,\\
            r_{k+1}&\leq r_k/(1-k_\zeta).
        \end{split}
    \end{align}
    Using \eqref{eq:d_and_r_inequalities} and noticing that $\gamma$ satisfies $\gamma/(1-k_\zeta)-\tau_c=(1-\tau_c)\gamma$, we obtain
    \begin{align*}
        s_{k+1} &\leq \frac{\gamma}{1-k_\zeta} r_k -(1-\tau_c)d_k -\tau_c r_k + \tau_c |\varrho_k| \\
        &\leq(1-\tau_c)(\gamma r_k-d_k)+ \tau_c |\varrho_k| \\
        \overset{\eqref{eq:def_aux_s_k}}&{=} (1-\tau_c) s_k + \tau_c |\varrho_k| \\
        \implies s_k &\leq (1-\tau_c)^{k} s_0 + \tau_c \sum_{i=0}^{k-1}(1-\tau_c)^{k-1-i}|\varrho_i|. \numberthis \label{eq:aux_s_k_bound}
    \end{align*}
    To obtain the eventual upper bound on $r_k$, we use \eqref{eq:aux_s_k_bound} to bound $\hat{d}_k$ from below. By \emph{Lem.~\ref{lemma:PE_DT}} and \eqref{eq:x_tilde_is_non_increasing}, $\varrho_k\to0$ exponentially fast. Define
    \begin{equation*}
        \scalemath{0.95}{\varepsilon_k :=
        (1-\tau_c)^k\max\{s_0,0\}
        +\tau_c\sum_{i=0}^{k-1}
        (1-\tau_c)^{k-1-i}|\varrho_i|
        +|\varrho_k|.}
    \end{equation*}
    Then $\varepsilon_k\to0$ and $\hat{d}_k\geq\gamma r_k-\varepsilon_k$. The reciprocal radius $\nu_k:=1/r_k$ then satisfies the linear comparison inequality
    \begin{align*}
        \nu_{k+1} \overset{\eqref{eq:reference_radius_update_DT}}&{=} \left(1-k_\zeta + \frac{k_\zeta}{d^*} \hat{d}_k\right) \nu_k \\
        \overset{\hat{d}_k\geq \gamma r_k - \varepsilon_k}&{\geq} \left(1-k_\zeta - \frac{k_\zeta \varepsilon_k }{d^*} \right) \nu_k + \frac{k_\zeta \gamma}{d^*} =: \varsigma_k \nu_k + \frac{k_\zeta \gamma}{d^*}.
    \end{align*}
    Since $\varepsilon_k\to 0$ as $k\to\infty$, for any fixed $\epsilon\in(0,1-k_\zeta)$, there exists some finite $\kappa_\epsilon \geq 0$ such that $\varsigma_k \geq 1-k_\zeta-\epsilon$ for all $k\geq \kappa_{\epsilon}$. Scalar linear comparison then gives
    \begin{equation*}
        \liminf_{k\to\infty}\nu_k \geq\frac{k_\zeta \gamma}{(k_\zeta+\epsilon)d^*}.
    \end{equation*}
    Letting $\epsilon\to 0^+$, we get
    \begin{equation}
        \label{eq:adaptive_reference_radius_bound_for_large_K}
        \scalemath{0.95}{\limsup_{k\to\infty}r_k \leq \frac{d^*}{\gamma} \overset{\eqref{eq:radial_compensator_gain_DT}}{<} \frac{3d^*}{2} \; \implies\; \exists \, K\geq 0 : r_k\leq \frac{3d^*}{2}, \; \forall\, k\geq K.}
    \end{equation}

    \textbf{Part II.} Introduce the transformed state $\bm{w}_k$,
    \begin{equation}
        \label{eq:coordinate_transform}
        \bm{w}_k :=\begin{bmatrix}
            p_k \\
            z_k
        \end{bmatrix}, \quad p_k:= \frac{d_k-d^*}{2d^*},\quad
        z_k:= \ln\frac{r_k}{r_{\infty}}.
    \end{equation}
    \emph{Step 1: SDC representation and coefficient bounds.} Let $\bar{d}_{k+1}$ and $\bar{r}_{k+1}$ denote the nominal one-step values obtained by setting $\varrho_k=0$ in \eqref{eq:def_d_k+1_integral} and \eqref{eq:reference_radius_update_DT}, evaluated at the current state $(d_k,r_k)$:
    \begin{subequations}
        \label{eq:nominal_maps}
        \begin{align}
            \bar{d}_{k+1}
            &:=
            \sqrt{\bigl((1-\tau_c)d_k+\tau_c r_k\bigr)^2
            +\tau_\omega^2},
            \label{eq:nominal_distance_map}\\
            \bar{r}_{k+1}
            &:=
            \frac{r_k}{1-k_\zeta+(k_\zeta/d^*)d_k}.
            \label{eq:nominal_compensator_update}
        \end{align}
    \end{subequations}
    Under \eqref{eq:coordinate_transform}, the updates \eqref{eq:nominal_maps} define the nominal one-step map $\bm{F}:\mathcal{D}\to\mathcal{D}$, $\mathcal{D}:=(-1/2,\infty)\times\mathbb{R}$, satisfying
    \begin{equation*}
        \bm{F}(p_k,z_k)=
        \begin{bmatrix}
            (\bar{d}_{k+1}-d^*)/(2d^*)\\[3pt]
            \ln(\bar{r}_{k+1}/r_\infty)
        \end{bmatrix},
    \end{equation*}
    and, expressed explicitly in terms of $(p,z)\in\mathcal{D}$,
    \begin{equation*}
        \scalemath{0.975}{
        \bm{F}(p,z)=
        \begin{bmatrix}
            \dfrac{
            \sqrt{\bigl((1-\tau_c)(1+2p)d^*
            +\tau_c r_\infty\me^z\bigr)^2+\tau_\omega^2}
            -d^*}{2d^*}\\[8pt]
            z-\ln(1+2k_\zeta p)
        \end{bmatrix}.}
    \end{equation*} 
    Observe that the map $\bm{F}\in\mathcal{C}^{1}$ on the convex domain $\mathcal{D}$ and satisfies $\bm{F}(0,0)=\bm{0}_2$. Therefore, the standard state-dependent coefficient (SDC) factorization, also known as \emph{extended linearization} (see, for example, \cite[Sec.~3.2]{cimen2010Systematic} or \cite[Sec.~3.1.3]{tsukamoto2021Contraction}), gives the following linear-like representation:
    \begin{equation}
        \label{eq:SDC}
        \bm{F}(p_k,z_k)=M_k\bm{w}_k,\quad
        M_k:=\int_0^1
        \mathrm{D}\bm{F}(\vartheta p_k, \,\vartheta z_k)\; \md \vartheta,
    \end{equation}
    where $\mathrm{D}\bm{F}$ denotes the Jacobian of $\bm{F}$. Let the difference between the actual and nominal transformed updates be $\bm{\rho}_k:=\bm{w}_{k+1}-\bm{F}(p_k,z_k)$, then the actual recursion becomes
    \begin{equation}
        \label{eq:radial_error_matrix_DT}
        \bm{w}_{k+1}=M_k\bm{w}_k+\bm{\rho}_k,
        \quad
        M_k=
        \begin{bmatrix}
            a_k & b_k\\
            -c_k & 1
        \end{bmatrix}.
    \end{equation} 
    Specifically, with
    \begin{equation*}
        \sigma_k(\vartheta)=(1-\tau_c)(1+2 p_{k} \vartheta)d^*+ \tau_c r_{\infty} \me^{ \vartheta z_{k}}, \quad 0 \leq \vartheta\leq 1,
    \end{equation*}
    the coefficients $a_k$, $b_k$, and $c_k$ are
    \begin{subequations}
        \label{eq:radial_error_matrix_bounds_DT}
        \begin{align}
            a_k &= (1-\tau_c)\int_0^1
               \frac{\sigma_k(\vartheta)}{\sqrt{(\sigma_k(\vartheta))^2+\tau_\omega^2}}\; \md \vartheta,\\
            b_k &= \frac{\tau_c r_\infty}{2d^*}\int_0^1
               \frac{\me^{ \vartheta z_k}\sigma_k(\vartheta)}
                    {\sqrt{(\sigma_k(\vartheta))^2+\tau_\omega^2}}\; \md \vartheta,\\
            c_k &= 2k_\zeta \int_0^1\frac{1}{1+2k_\zeta \vartheta p_k}\; \md \vartheta.
        \end{align}
    \end{subequations}
    Here $p_k>-1/2$, and both $\sigma_k(\vartheta)$ and $1+2k_\zeta \vartheta p_k$ are positive for $\vartheta\in[0,1]$. Set $S_k=\sigma_k(1)=(1-\tau_c)d_k+ \tau_c r_k$. The perturbation is
    \begin{equation}
        \label{eq:def_rho_k}
        \bm{\rho}_k =
        \begin{bmatrix}
            \dfrac{\sqrt{(S_k+\tau_c \varrho_k)^2+\tau_\omega^2}
              -\sqrt{S_k^2+\tau_\omega^2}}{2d^*}\\[10pt]
            \ln\dfrac{1-k_\zeta+(k_\zeta/d^*)d_k}{1-k_\zeta+(k_\zeta/d^*)(d_k-\varrho_k)}
        \end{bmatrix}.
    \end{equation}
    Since $d_k\geq0$ and $d_k-\varrho_k=\hat d_k\geq0$, the numerator and denominator inside the logarithm are both bounded below by $1-k_\zeta>0$. Applying the triangle inequality and the mean value theorem gives
    \begin{equation}
        \label{eq:rho_k_bound}
        \scalemath{0.95}{\|\bm{\rho}_k\| \leq \varpi |\varrho_k|, \quad \varpi:= \sqrt{ 
        \left( \frac{\tau_c}{2d^*} \right)^2 + \left( \frac{k_\zeta}{d^*(1-k_\zeta)} \right)^2
        } \; .}
    \end{equation}
For $k\geq K$, Part I and $r_\infty<d^*$ imply $r_\infty \me^{\vartheta z_k}=(r_{\infty})^{1-\vartheta}r_k^{\vartheta}\leq 3d^*/2$ for $\vartheta\in [0,1]$. Hence the coefficients satisfy the following bound,
\begin{subequations}
    \label{eq:corrected_tail_matrix_bounds_DT}
    \begin{gather}
        \begin{split}
            0\leq a_k\leq1-\tau_c ,\quad
            0<b_k\leq\frac{3\tau_c}{4},
        \end{split} \\
        0 < c_k \leq \bar{c}:=\frac{2k_\zeta}{1-k_\zeta} \overset{\eqref{eq:radial_compensator_gain_DT}}{<} \tau_c.
    \end{gather}
\end{subequations}

\begin{figure}[t]
        \centering
        \includegraphics[width=0.865\linewidth]{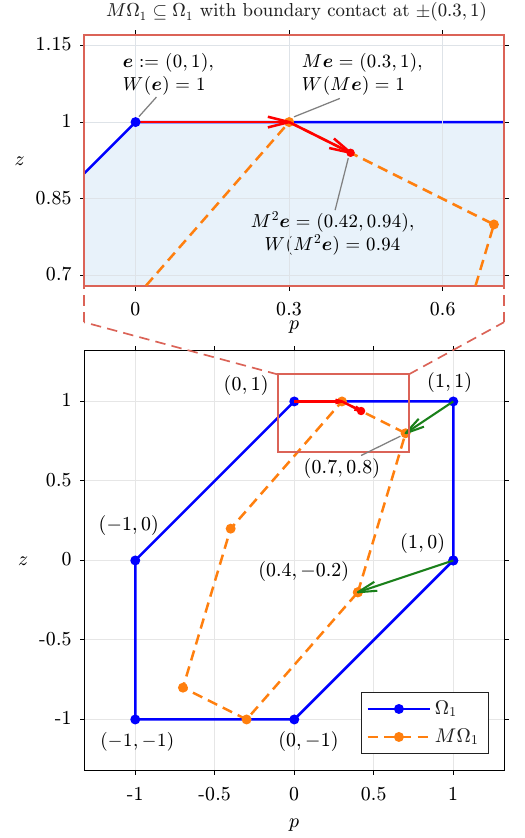}
        \caption{\label{fig:Omega_illustration}The unit sublevel set $\Omega_1$ and its image $M\Omega_1:=\{M\bm{\nu}:\,\bm{\nu}\in\Omega_1\}$ under the linear transformation $\bm{\nu}\mapsto M\bm{\nu}$. Here,
        $M:=\begin{bmatrix*}[r]
            0.4 & 0.3 \\
            -0.2 & 1
        \end{bmatrix*}$
        is a constant matrix, held fixed over the two-step evolution. Although the numerical entries of $M$ are chosen for illustrative purposes, $M$ has the same structure as $M_k$ in \eqref{eq:radial_error_matrix_DT} and satisfies the coefficient bounds in \eqref{eq:corrected_tail_matrix_bounds_DT} for the illustrative parameter values $\tau_c=0.5$ and $k_\zeta=0.1$.}
    \end{figure}

\noindent \emph{Step 2: Non-expansiveness and boundedness.} Consider the polyhedral Lyapunov function (see, for instance, \cite[Sec.~4.5]{blanchini2015SetTheoretic})
\begin{equation}
    \label{eq:polyhedral_Lyapunov}
    W(p,z)=\max\{|p|,|z|,|p-z|\},
\end{equation}
whose unit sublevel set $\Omega_1:=\{\bm{\nu}\in\mathbb{R}^2: \, W(\bm{\nu})\leq 1\}$ is a hexagon and has vertices $\pm[1,0]^\top$, $\pm[1,1]^\top$, and $\pm[0,1]^\top$. To aid the discussion, a graphical illustration of $\Omega_1$ is provided in Fig.~\ref{fig:Omega_illustration}. For $k\geq K$, \eqref{eq:corrected_tail_matrix_bounds_DT} applies and the images of these vertices under $M_k$ satisfy
\begin{align*}
    W\bigl(M_k[1,0]^\top \bigr) &=a_k+c_k\leq 1-\tau_c +\bar{c}<1,\\
    W\bigl(M_k[1,1]^\top \bigr) &= \max\{a_k+b_k,1-c_k\}\\
    &\leq\max\{1-\tau_c/4,1-c_k\}<1,\\
    W\bigl(M_k[0,1]^\top \bigr) &=1.
\end{align*}
By symmetry, convexity, and homogeneity, $W(M_k\bm{\nu})\leq W(\bm{\nu})$ for every $\bm{\nu}\in\mathbb R^2$ and $k\geq K$. Using \eqref{eq:radial_error_matrix_DT}, the bound on $\bm{\rho}_k$ \eqref{eq:rho_k_bound}, and $W(\bm{\nu})\leq\sqrt{2}\|\bm{\nu}\|$, we obtain
\begin{equation*}
    W(\bm w_{k+1})\leq W(\bm w_k)+\sqrt2\varpi|\varrho_k|,\quad k\geq K.
\end{equation*}
Since $\sum_{k=0}^{\infty} \|\bm{\rho}_k\| \leq \varpi \sum_{k=0}^{\infty} |\varrho_k| <\infty$, the sequence $W(\bm{w}_k)$ is bounded. In particular, $r_k= r_\infty \me^{z_k}$ is bounded away from zero. Together with the preceding bounds on $d_k$ and $r_k$, this places $(d_k,r_k)$ in a compact subset of $(0,\infty)\times (0,\infty)$ for all subsequent $k\geq K$. Continuity and strict positivity of the functions defining $b_k$ and $c_k$ in \eqref{eq:radial_error_matrix_bounds_DT} now imply that there exist constants $\underline{b},\underline{c}>0$ such that
\begin{equation}
    \label{eq:k_uniform_lower_bounds}
    b_k \geq \underline{b}, \quad c_k \geq \underline{c}, \quad \forall\, k \geq K.
\end{equation}
 
    \noindent \emph{Step 3: Two-step contraction and exponential convergence.} The only non-contractive vertices are $[0,\pm 1]^\top$. For $[0,1]^\top$, two consecutive steps give
    \begin{equation}
        \label{eq:consecutive_steps}
        M_{k+1}M_k
        \begin{bmatrix}0\\1\end{bmatrix}
        =
        \begin{bmatrix}
            a_{k+1}b_k+b_{k+1}\\
            1-c_{k+1}b_k
        \end{bmatrix}.
    \end{equation}
    By \eqref{eq:corrected_tail_matrix_bounds_DT} and \eqref{eq:k_uniform_lower_bounds}, the two entries in \eqref{eq:consecutive_steps} satisfy the following inequalities, uniformly over all $k\geq K$,
    \begin{align*}
        0
        &\leq a_{k+1}b_k+b_{k+1}
        \leq \frac{3\tau_c}{4}(2-\tau_c)
        \leq \frac{3}{4},\\
        0
        &< 1-c_{k+1}b_k
        \leq 1-\underline c\,\underline b
        < 1.
    \end{align*}
    By construction \eqref{eq:polyhedral_Lyapunov}, $W(p,z)=\max\{p,z\}$ for $p,z\geq0$, it then follows, for all $k\geq K$, that
    \begin{align*}
        W\bigl(M_{k+1}M_k[0,1]^\top\bigr)
        &= \max\bigl\{
            a_{k+1}b_k+b_{k+1},\,
            1-c_{k+1}b_k
           \bigr\}\\
        &\leq
        \max\left\{
            \frac{3}{4},\;
            1-\underline c\,\underline b
        \right\}
        <1.
    \end{align*}
    Observe that the upper bound is independent of $k$. A graphical illustration involving a sample two-step evolution from the vertex $[0,1]^\top$ is shown using red arrows in Fig.~\ref{fig:Omega_illustration} to aid the discussion. By symmetry, the same bound holds for $W\bigl(M_{k+1}M_k[0,-1]^\top\bigr)$. The same vertex argument therefore implies the existence of a constant $\theta\in(0,1)$ such that
    \begin{equation*}
        W(M_{k+1}M_k\bm{\nu})\leq\theta W(\bm{\nu}), \quad \forall\, k\geq K,\; \forall\, \bm{\nu}\in\mathbb{R}^2.
    \end{equation*}
    Applying \eqref{eq:radial_error_matrix_DT} twice gives
    \begin{equation*}
        \bm w_{k+2}=M_{k+1}M_k\bm w_k+M_{k+1}\bm\rho_k+\bm\rho_{k+1},
    \end{equation*}
    and therefore, we obtain
    \begin{equation*}
        W(\bm w_{k+2})\leq\theta W(\bm w_k)
        +\sqrt2\varpi\bigl(|\varrho_k|+|\varrho_{k+1}|\bigr), \quad \forall \, k\geq K.
    \end{equation*}
    Exponential decay of $\varrho_k$ now implies exponential convergence of $p_k,z_k$ to zero. Since $d_k-d^*=2d^* p_k$ and $r_k-r_\infty =r_{\infty} (\me^{z_k}-1)$, the claimed convergence follows.
\end{proof}

    \begin{figure*}[ht]
        \centering
        \includegraphics[width=1.0\linewidth]{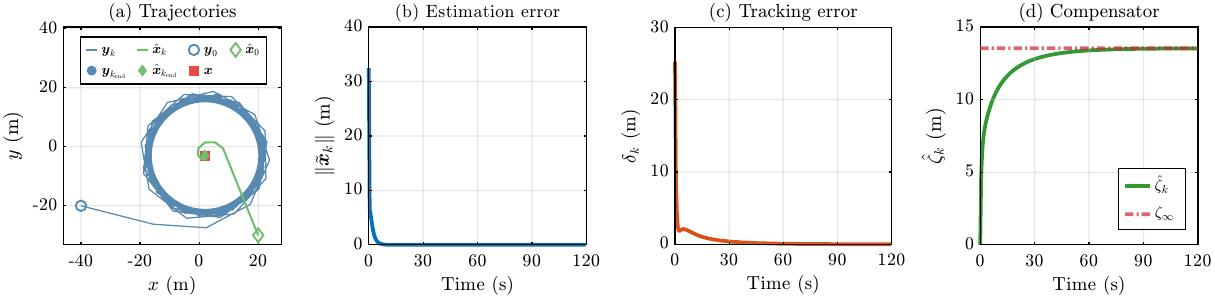}
        \caption{\label{fig:stationary_target_sim_plot}Simulation results of the proposed circumnavigation controller with integral action \eqref{eq:def_controller_integral_DT} coupled with the target estimator \eqref{eq:def_target_estimator_DT}.}
    \end{figure*}

\section{Simulations}
\label{sec:Simulations}

    We present the following simulation study. The evolution of the system \eqref{eq:agent_kinematics_DT}, driven by the estimator \eqref{eq:def_target_estimator_DT} and the proposed controller with integral action \eqref{eq:def_controller_integral_DT}, is simulated in MATLAB by direct iteration over $120$~s with sampling interval $T=0.5$~s (240 control updates). The setup is: $\bm{x}_k\equiv[2,\,-3]^\top\,\mathrm{m}$, $\bm{y}_0=[-40,\,-20]^\top\,\mathrm{m}$, $\hat{\bm{x}}_0=[20,\,-30]^\top\,\mathrm{m}$, $d^*=20\,\mathrm{m}$, $\hat{\zeta}_0=0$, $k_c=k_x=1$, $k_\omega=30$, and $k_\zeta=0.15$. Figure~\ref{fig:stationary_target_sim_plot}(a) presents the resulting agent trajectory and target estimates. Figures~\ref{fig:stationary_target_sim_plot}(b) and \ref{fig:stationary_target_sim_plot}(c) show that the estimation error $\|\tilde{\bm{x}}_k\|$ and sampled tracking error $\delta_k$ approach zero, which are consistent with \emph{Lemma~\ref{lemma:PE_DT}} and \emph{Theorem~\ref{theorem:circum_error_goes_to_zero}}. Fig.~\ref{fig:stationary_target_sim_plot}(d) shows that the adaptive compensator $\hat{\zeta}_k$ approaches the predicted limit $\zeta_\infty\approx 13.542\,\mathrm{m}$ given by \eqref{eq:def_zeta_ss_geometric}, thereby illustrating that the proposed compensator generates the required inward correction without explicit knowledge of $T$.

\section{CONCLUDING REMARKS}
\label{sec:conclusion}
    In this paper, an adaptive compensator and an improved discrete-time controller were proposed for the localization and circumnavigation of an unknown target using bearing measurements and the agent's own position. With the help of the adaptive compensator, the estimator-coupled controller attained zero steady-state error in radial tracking at the sampling instants. The required geometric correction was learned dynamically without explicitly using the sampling interval in either the control law or the compensator update. Exponential convergence of the target estimation and sampled tracking errors was rigorously established under the stated sampling and gain conditions. Simulation results illustrate satisfactory performance of the proposed algorithms.

    Future directions include extending the proposed method to moving targets, accounting for measurement noise and non-uniform sampling, and evaluating the method in field experiments.

    \noindent \emph{Acknowledgment.} The authors would like to thank Shuo Xu and Zhiyong Sun for their helpful comments on the early version of the manuscript and valuable discussions.


\addcontentsline{toc}{chapter}{Bibliography}
\bibliographystyle{ieeetr}
\bibliography{myRefs}

\end{document}